\documentclass[letterpaper, 10pt, conference]{ieeeconf}  

\IEEEoverridecommandlockouts           

\usepackage{amsmath} 	
\usepackage{amssymb}  	
\usepackage{amsfonts}

\usepackage{graphicx}
\usepackage{subfigure}
\usepackage{epstopdf}
\usepackage{dblfloatfix}
\usepackage{balance}

\usepackage{cite}

\usepackage{tikz}
\usetikzlibrary{shapes,matrix,decorations.shapes}

\usepackage{array,xcolor}

\usepackage{algorithm}
\usepackage{algpseudocode}

\usepackage{multirow}
\usepackage{threeparttable}
\usepackage{booktabs}

\usepackage[colorlinks=true,      
			linkcolor=black,      
			citecolor=black,      
			filecolor=black,      
			urlcolor=blue ]{hyperref}

\newtheorem{definition}{Definition}
\newtheorem{theorem}{Theorem}

\newtheorem{remark}{Remark}

\title{\LARGE \bf
Fast ADMM for Semidefinite Programs with Chordal Sparsity
}

\author{Yang~Zheng$^{\dagger,1}$,~Giovanni~Fantuzzi$^{\dagger,2}$,~Antonis Papachristodoulou$^{1}$,~Paul Goulart$^{1}$,~Andrew~Wynn$^{2}$ 
\thanks{$^\dagger$Y.~Zheng and G.~Fantuzzi contributed equally to this work. Y. Zheng is supported by the Clarendon Scholarship and the Jason Hu Scholarship. G.~Fantuzzi was partially supported by the EPSRC grant EP/J010537/1.}
\thanks{$^1$Department of Engineering Science, University of Oxford, Parks Road, Oxford, OX1 3PJ, United Kingdom (e-mail:
yang.zheng@eng.ox.ac.uk; paul.goulart@eng.ox.ac.uk; antonis@eng.ox.ac.uk).}%
\thanks{$^2$Department of Aeronautics, Imperial College London, South Kensington Campus, London, SW7 2AZ, United Kingdom (e-mail: gf910@ic.ac.uk; a.wynn@imperial.ac.uk).}%
}

\begin{document}

\maketitle
\thispagestyle{empty}
\pagestyle{empty}

\begin{abstract}
	Many problems in control theory can be formulated as semidefinite programs (SDPs). For large-scale SDPs, it is important to exploit the inherent sparsity to improve the scalability. This paper develops efficient first-order methods to solve SDPs with chordal sparsity based on the alternating direction method of multipliers (ADMM). We show that chordal decomposition can be applied to either the primal or the dual standard form of a sparse SDP, resulting in scaled versions of ADMM algorithms with the same computational cost. Each iteration of our algorithms consists of a projection on the product of small positive semidefinite cones, followed by a projection on an affine set, both of which can be carried out efficiently. Our techniques are implemented in CDCS, an open source add-on to MATLAB. Numerical experiments on large-scale sparse problems in SDPLIB and random SDPs with block-arrow sparse patterns show speedups compared to some common state-of-the-art software packages.
\end{abstract}

\section{Introduction}

Semidefinite programs (SDPs) are a type of convex optimization problems over the cone of positive semidefinite (PSD) matrices.
Given $ b\in \mathbb{R}^m$, $C\in \mathbb{S}^n$, and matrices $A_1,\,\ldots,\,A_m \in \mathbb{S}^n$ that define the operators
\begin{equation*}
\mathcal{A}(X) =
\begin{bmatrix}
\langle A_1,X \rangle\\ \vdots \\ \langle A_m,X \rangle
\end{bmatrix},
\quad
\mathcal{A}^*(y) = \sum_{i=1}^m A_i y_i,
\end{equation*}
SDPs are typically written in the \textit{standard primal form}
\begin{equation}
\label{E:PrimalSDP}
    \begin{aligned}
    \min_{X} \quad & \langle C, X \rangle \\
    \text{subject to } \quad & \mathcal{A}(X) = b,\\
        & X \in \mathbb{S}^n_{+},
    \end{aligned}
\end{equation}
or in the \textit{standard dual form}
\begin{equation}
\label{E:DualSDP}
    \begin{aligned}
        \max_{y,Z} \quad & \langle b,y \rangle\\
        \text{subject to } \quad  & \mathcal{A}^*(y) + Z= C,\\
        & Z \in \mathbb{S}^n_{+}.
    \end{aligned}\\
\end{equation}
In the above and throughout this work, $\mathbb{R}^m$ is the $m$-dimensional Euclidean space, $\mathbb{S}^n$ is the space of $n\!\times\!n$ symmetric matrices, $\mathbb{S}^n_{+}$ is the subspace of PSD matrices,
and $\langle \cdot,\cdot \rangle$ denotes the inner product in the appropriate space.

SDPs have applications in control theory, machine learning, combinatorics, and operations research~\cite{boyd1994linear}. Moreover, linear, quadratic, and second-order-cone programs, are particular instances of SDPs~\cite{boyd2004convex}.
Small to medium-sized SDPs can be solved in polynomial time~\cite{vandenberghe1996semidefinite} using efficient second-order interior-point methods (IPMs)~\cite{alizadeh1998primal, helmberg1996interior}. However, many real-life problems are too large for the state-of-the-art interior-point algorithms, due to memory and CPU time constraints.

One approach is to abandon IPMs, in favour of faster first-order methods (FOMs) with modest accuracy. For instance, Wen \emph{et al.} proposed an alternating-direction augmented Lagrangian method for large-scale SDPs in the dual standard form~\cite{wen2010alternating}. More recently, O'Donoghue \emph{et al.} developed a first-order operator-splitting method to solve the homogeneous self-dual embedding (HSDE) of a primal-dual pair of conic programs, which has the advantage of being able to provide primal or dual certificates of infeasibility~\cite{ODonoghue2016}. 
A second approach relies on the fact that large-scale SDPs are often structured and/or sparse~\cite{boyd1994linear}. Exploiting sparsity in SDPs is an active and challenging area of research~\cite{andersen2011interior}, one main difficulty being that the optimal solution is typically dense despite the sparsity of the problem data. If, however, the sparsity pattern of the data is \textit{chordal} or has sparse \textit{chordal extensions}, 
Grone's and Agler's theorems~\cite{grone1984positive,agler1988positive} allow replacing the PSD constraint with a set of smaller semidefinite constraints, plus an additional set of equality constraints. In some cases, the converted SDP can then be solved more efficiently than the original problem. These ideas underly the \emph{domain-} and \emph{range-space} conversion techniques~\cite{fukuda2001exploiting,kim2011exploiting}, implemented in SparseCoLO~\cite{fujisawa2009user}.

However, adding equality constraints often offsets the benefit of working with smaller PSD cones. One possible solution is to exploit chordal sparsity directly in the IPMs: Fukuda \emph{et al.} used Grone's theorem~\cite{grone1984positive} to develop a primal-dual path-following method for SDPs~\cite{fukuda2001exploiting}; Burer proposed a nonsymmetric primal-dual IPM using Cholesky factors of the dual variable and maximum determinant completion of the primal variable~\cite{burer2003semidefinite}; and Andersen \emph{et al.} developed fast recursive algorithms for SDPs with chordal sparsity~\cite{andersen2010implementation}. Alternatively, one can solve the decomposed SDP with FOMs: Sun \emph{et al.} proposed a first-order splitting method for decomposable conic programs~\cite{sun2014decomposition};
Kalbat \& Lavaei applied the alternating-direction method of multipliers (ADMM) to SDPs with fully decomposable constraints~\cite{Kalbat2015Fast}; Madani \emph{et al.} developed a highly-parallelizable ADMM algorithm for sparse SDPs with inequality constraints with optimal power flow applications~\cite{Madani2015ADMM}.

In this work we adopt the strategy of exploiting sparsity using first-order algorithms in the spirit of~\cite{sun2014decomposition, Kalbat2015Fast, Madani2015ADMM}, and develop efficient ADMM algorithms to solve large-scale sparse SDPs. 
%
Our contributions are:
\begin{enumerate}
  \item We combine ADMM and chordal decomposition to solve sparse SDPs in primal or dual standard form. The resulting algorithms are scaled versions of each other. This gives a conversion framework for the application of FOMs, analogous to that of~\cite{fukuda2001exploiting, kim2011exploiting} for IPMs. 
  \item In each iteration, the PSD constraint
is enforced via parallel projections onto small PSD cones. The affine constraints
are imposed by a quadratic program with equality constraints, and its KKT system matrix can be factorized before iterating the ADMM algorithm since it only depends on the problem data.
  \item We implement our methods in the open-source MATLAB solver CDCS (Cone Decomposition Conic Solver)~\cite{CDCS}.  Numerical simulations on random SDPs with block-arrow sparsity patterns and on four large-scale sparse problems in SDPLIB \cite{borchers1999sdplib} demonstrate the efficiency of our algorithms compared to other solvers.
\end{enumerate}

The rest of this paper is organized as follows. Section~\ref{se:preliminaries} reviews chordal sparsity and decomposition techniques. We show how to apply the ADMM to primal and dual standard-form SDPs in Sections~\ref{se:pSDP}--\ref{se:dSDP}, respectively, and report our numerical experiments in Section~\ref{se:simulation}. Finally, Section~\ref{se:conclusion} offers concluding remarks.

\section{Preliminaries: Chordal Decomposition and the ADMM Algorithm} \label{se:preliminaries}

\subsection{Chordal graphs}

Let $\mathcal{G}(\mathcal{V},\mathcal{E})$ be an undirected graph with vertices $\mathcal{V}=\{1,2,\ldots,n\}$ and edges $\mathcal{E} \subseteq \mathcal{V} \times \mathcal{V}$.
A \emph{clique} $\mathcal{C}\subseteq \mathcal{V} $ is a subset of vertices such that $ (i,j) \in \mathcal{E}$ for any distinct vertices $ i,j \in \mathcal{C} $, and the number of vertices in $\mathcal{C}$ is denoted by $\vert \mathcal{C} \vert$. If $\mathcal{C}$ is not a subset of any other clique, then it is referred to as a \emph{maximal clique}.
A cycle of length $ k $ in $\mathcal{G}$ is a set of pairwise distinct vertices $ \{v_1,v_2,\ldots,v_k\}\subset\mathcal{V} $ such that $ (v_k,v_1) \in \mathcal{E} $ and $ (v_i,v_{i+1}) \in \mathcal{E} $ for $ i=1,\ldots,k-1 $. A chord is an edge joining two non-adjacent vertices in a cycle.

\begin{definition}[Chordal graph]
An undirected graph is \emph{chordal} if all cycles of length four or higher have a chord.
\end{definition}

Note that if $\mathcal{G}(\mathcal{V}, \mathcal{E})$ is not chordal, it can be \emph{chordal extended}, \emph{i.e.}, we can construct a chordal graph $\mathcal{G}'(\mathcal{V}, \mathcal{E}') $ by adding edges to $ \mathcal{E} $ such that $\mathcal{G}'$ is chordal. Finding the chordal extension with the minimum number of additional edges is an NP-complete problem~\cite{yannakakis1981computing}, but good chordal extensions can be computed efficiently using several heuristics~\cite{vandenberghe2014chordal}.

\subsection{Sparse matrices defined by graphs}
Let $ \mathcal{G} = (\mathcal{V},\mathcal{E}) $ be an undirected graph such that $(i,i) \in \mathcal{E}$, \emph{i.e.}, each node has a self-loop. We say that $X$ is a sparse symmetric matrix defined by $\mathcal{G}$ if $X_{ij}=X_{ji} = 0$ whenever $(i,j)\notin\mathcal{E}$. The spaces of sparse and PSD sparse symmetric matrices defined by $\mathcal{G}$ are
\begin{equation*}
\begin{aligned}
\mathbb{S}^n(\mathcal{E},0) = &\{ X \in \mathbb{S}^n \mid X_{ij} =X_{ji} = 0 \text{ if } (i,j) \notin \mathcal{E}  \}, \\
\mathbb{S}_{+}^n(\mathcal{E},0) = &\{ X \in \mathbb{S}^n(\mathcal{E},0) \mid X \succeq 0 \}.
\end{aligned}
\end{equation*}
Similarly, we say that $X$ is a partial symmetric matrix defined by $\mathcal{G}$ if $X_{ij}=X_{ji}$ are given when $(i,j)\in\mathcal{E}$, and arbitrary otherwise. Moreover, we say that $M$ is a PSD completion of the partial symmetric matrix $X$ if  $M\succeq 0$ and $M_{ij}=X_{ij}$ when $(i,j)\in\mathcal{E}$. We can then define the spaces
\begin{equation*}
\begin{aligned}
\mathbb{S}^n(\mathcal{E},?)\!= &\{ X\!\in \mathbb{S}^n \mid X_{ij} =X_{ji}\text{ given if } (i,j) \in \mathcal{E}  \}, \\
\mathbb{S}_{+}^n(\mathcal{E},?)\!= &\{ X\!\in\mathbb{S}^n(\mathcal{E},?) \mid
\!\exists M \succeq 0, \, M_{ij}\!=\!X_{ij}\, \forall (i,j)\!\in \mathcal{E}  \}.
\end{aligned}
\end{equation*}

Finally, given a clique $\mathcal{C}_k$ of $\mathcal{G}$, $E_k \in \mathbb{R}^{\mid \mathcal{C}_k\mid \times n}$ is the matrix with $(E_{k})_{ij} = 1$ if $\mathcal{C}_k(i) = j$ and zero otherwise, where $\mathcal{C}_k(i)$  is the $i$-th vertex in $\mathcal{C}_k$, sorted in the natural ordering. The submatrix of $X \in \mathbb{S}^n$ defined by $\mathcal{C}_k$ is $E_{k}XE_{k}^T \in \mathbb{S}^{\mid \mathcal{C}_k\mid}$.

\subsection{Chordal decomposition of PSD matrices}

The spaces $\mathbb{S}^n_{+}(\mathcal{E},?)$ and $\mathbb{S}_{+}^n(\mathcal{E},0)$ are a pair of dual convex cones for any undirected graph $\mathcal{G}(\mathcal{V},\mathcal{E})$~\cite{andersen2010implementation,vandenberghe2014chordal}. If $\mathcal{G}$ is chordal, $\mathbb{S}^n_{+}(\mathcal{E},?)$ and $\mathbb{S}_{+}^n(\mathcal{E},0)$ can be expressed in terms of several coupled smaller convex cones:
\begin{theorem} [Grone's theorem~\cite{grone1984positive}]\label{T:ChordalCompletionTheorem}
     Let $\{\mathcal{C}_1,\mathcal{C}_2, \ldots, \mathcal{C}_p\}$  be the set of maximal cliques of a chordal graph $\mathcal{G}(\mathcal{V},\mathcal{E})$. Then, $X\in\mathbb{S}^n_+(\mathcal{E},?)$ if and only if
    $X_k := E_k X E_k^T \in \mathbb{S}^{\vert \mathcal{C}_k \vert}_+$
    for all $k=1,\,\ldots,\,p$.
\end{theorem}
\begin{theorem} [Agler's theorem~\cite{agler1988positive}]\label{T:ChordalDecompositionTheorem}
     Let $\{\mathcal{C}_1,\mathcal{C}_2, \ldots, \mathcal{C}_p\}$  be the set of maximal cliques of a chordal graph $\mathcal{G}(\mathcal{V},\mathcal{E})$. Then, $Z\in\mathbb{S}^n_+(\mathcal{E},0)$ if and only if there exist matrices $Z_k \in \mathbb{S}^{\vert \mathcal{C}_k \vert}_+$ for $k=1,\,\ldots,\,p$ such that
    $Z = \sum_{k=1}^{p} E_k^T Z_k E_k.$
\end{theorem}

These results can be proven individually, but can also can be derived from each other using the duality of the cones $\mathbb{S}^n_{+}(\mathcal{E},?)$ and $\mathbb{S}_{+}^n(\mathcal{E},0)$~\cite{vandenberghe2014chordal}.



\subsection{ADMM algorithm}

The ADMM algorithm solves the optimization problem
\begin{equation*}
    \begin{aligned}
        \min \quad & f(x)+g(y) \\
        \text{subject to } \quad & Ax + By = c,
    \end{aligned}
\end{equation*}
where $f$ and $g$ are convex functions, $x \in \mathbb{R}^{n_x}, y \in \mathbb{R}^{n_y}, A \in \mathbb{R}^{n_c\times n_x}, B \in \mathbb{R}^{n_c\times n_y}$ and $c \in \mathbb{R}^{n_c}$. Given a penalty parameter $\rho>0$ and a dual multiplier $z \in \mathbb{R}^{n_c}$, the ADMM algorithm minimizes the augmented Lagrangian
\begin{equation*}
 L_{\rho}(x,y,z) = f(x) + g(y)
 + \frac{\rho}{2} \left\|Ax + By - c + \frac{1}{\rho} z\right\|^2
\end{equation*}
with respect to the variables $x$ and $y$ separately, followed by a dual variable update:
\begin{subequations}\label{E:ADMM}
    \begin{align}
        x^{(n+1)} & = \text{arg} \min_{x} L_{\rho}(x,y^{(n)},z^{(n)}),
        \label{E:ADMM_S1}\\
        y^{(n+1)} & = \text{arg} \min_{y} L_{\rho}(x^{(n+1)},y,z^{(n)}),
        \label{E:ADMM_S2}\\
        z^{(n+1)} &= z^{(n)} + \rho ( A x^{(n+1)} + B y^{(n+1)} - c).
    \end{align}
\end{subequations}
The superscript $(n)$ indicates that a variable is fixed to its value at the $n$-th iteration. ADMM is particularly suitable when the minimizations with respect to each of the variables $x$ and $y$ in~\eqref{E:ADMM_S1} and~\eqref{E:ADMM_S2} can be carried out efficiently through closed-form expressions.

\section{ADMM for Sparse Primal-form SDPs} \label{se:pSDP}


\subsection{Reformulation and decomposition of the PSD constraint}

Let the primal-standard-form SDP~\eqref{E:PrimalSDP} be sparse with an \emph{aggregate sparsity pattern} described by the graph $\mathcal{G}(\mathcal{V},\mathcal{E})$, meaning that $(i,j)\in\mathcal{E}$ if and only if the entry $ij$ of at least one of the data matrices $C,\,A_0,\,\ldots,\,A_m$, is nonzero. We assume that $\mathcal{G}$ is chordal (otherwise it can be chordal extended) and that its maximal cliques $\mathcal{C}_1,\ldots,\mathcal{C}_p$ are small. Then, only the entries of the matrix variable $X$ corresponding to the graph edges $\mathcal{E}$ appear in the cost and constraint functions, so the constraint $X\in\mathbb{S}^n_+$ can be replaced by $X\in\mathbb{S}^n_+(\mathcal{E},?)$. Using Theorem~\ref{T:ChordalCompletionTheorem}, we can then reformulate~\eqref{E:PrimalSDP} as
\begin{equation}
\label{E:DecomposedPrimalSDP}
\begin{aligned}
    \min_{X,X_1,\ldots,X_p} \quad & \langle C, X \rangle \\[-0.25em]
    \text{subject to } \quad  &\mathcal{A}(X) = b,\\
		& X_k - E_k X E_k^T = 0, &&k=1,\,\ldots,\,p,\\
        & X_k \in \mathbb{S}^{\vert \mathcal{C}_k\vert}_{+}, &&k=1,\,\ldots,\,p.
\end{aligned}
\end{equation}
In other words, we can decompose the original large semidefinite cone into multiple smaller cones, at the expense of introducing a set of consensus constraints between the variables.

To ease the exposition, we rewrite~\eqref{E:DecomposedPrimalSDP} in a vectorized form. Letting $\mathrm{vec}:\mathbb{S}^n \to \mathbb{R}^{n^2}$ be the usual operator mapping a matrix to the stack of its column, define the vectorized data
%
$ 
    c := \text{vec}(C),
    A := \begin{bmatrix} \text{vec}(A_0) & \hdots & \text{vec}(A_m) \end{bmatrix}^T,
$
the vectorized variables
%
$
        x := \mathrm{vec}(X),    x_k := \text{vec}(x_k),   \quad k = 1,\,\ldots,\,p,
$
and the matrices $H_k := E_k \otimes E_k$ such that
$x_k = \mathrm{vec}(X_k) = \mathrm{vec}(E_k X E_k^T) = H_k x$.
%
%
In other words, the matrices $H_1,\,\ldots,\,H_p$ are ``entry-selector'' matrices of $1$'s  and $0$'s, whose rows are orthonormal, that project $x$ onto the subvectors $x_1,\,\ldots,\,x_p$, respectively.
Denoting the constraints $X_k\in\mathbb{S}^{|\mathcal{C}_k|}_+$ by $x_k\in\mathcal{S}_k$, we can rewrite~\eqref{E:DecomposedPrimalSDP} as
\begin{equation} \label{E:PrimalVectorForm}
\begin{aligned}
\min_{x,x_1,\ldots, x_p} \quad & \langle c,x \rangle\\[-0.25em]
\text{subject to } \quad & Ax=b,\\
					   & x_k = H_k x,&& k=1,\,\ldots,\,p,\\
					   & x_k \in \mathcal{S}_k,&& k=1,\,\ldots,\,p.
\end{aligned}
\end{equation}

\subsection{The ADMM algorithm for primal SDPs}
\label{S:ADMMAlgorithmPrimal}

Moving the constraints $Ax = b$ and $x_k \in \mathcal{S}_k$ in \eqref{E:PrimalVectorForm} to the objective using the indicator functions $\delta_0(\cdot)$ and $\delta_{\mathcal{S}_k}(\cdot)$ gives
\begin{equation} \label{E:ADMMPrimal}
    \begin{aligned}
    \min_{x,x_1,\ldots,x_p} \quad &\langle c,x \rangle + \delta_0\left( Ax - b \right) + \sum_{k=1}^{p} \delta_{\mathcal{S}_k}(x_k)
    \\
    \text{subject to } \quad &x_k = H_k x, \quad k=1,\,\ldots,\,p.
    \end{aligned}
\end{equation}

This problem is in the standard form for the application of ADMM. Given a penalty parameter $\rho>0$ and a Lagrange multiplier $\lambda_k$ for each constraint $x_k = H_k x$, we define
\begin{multline}
\label{E:AugLagrPrimal}
\mathcal{L}:= \langle c,x \rangle + \delta_0\left( Ax - b \right)
\\
+ \sum_{k=1}^{p} \left[ \delta_{\mathcal{S}_k}(x_k) + \frac{\rho}{2}\left\| x_k - H_k x + \frac{1}{\rho}\lambda_k \right\|^2 \right],
\end{multline}
and group the variables as
    $\mathcal{X}:= \{x\}$,
    $\mathcal{Y}:= \{x_1,\,\ldots,\,x_p\}$, and
    $\mathcal{Z}:= \{\lambda_1,\,\ldots,\,\lambda_p\}$.
As in~\eqref{E:ADMM}, in each ADMM iteration, we minimize~\eqref{E:AugLagrPrimal} with respect to $\mathcal{X}$ and $\mathcal{Y}$, then update $\mathcal{Z}$. 

\subsubsection{\textbf{Minimization over $\mathcal{X}$}} \label{se:MinYblkprimal}

Minimizing~\eqref{E:AugLagrPrimal} over $\mathcal{X}$ is equivalent to the equality-constrained quadratic program
\begin{equation} \label{E:MinXblockPrimal}
    \begin{aligned}
        \min_{x} \quad &\langle c,x \rangle + \frac{\rho}{2}\sum_{k=1}^{p} \left\| x_k^{(n)} - H_k x + \frac{1}{\rho}\lambda_k^{(n)} \right\|^2
        \\
        \text{subject to } \quad &Ax=b.
    \end{aligned}
\end{equation}
Define
$
D := \sum_{k=1}^{p} H_k^T H_k
$
and let $\rho y$ be the multiplier for the equality constraint. We can write the optimality conditions for~\eqref{E:MinXblockPrimal} as the KKT system
\begin{equation} \label{E:OptCondMinYPrimal}
    \begin{bmatrix}D & A^T \\ A & 0\end{bmatrix}
    \begin{bmatrix}x \\ y\end{bmatrix} =
    \begin{bmatrix}\sum_{k=1}^{p} H_k^T\left( x_k^{(n)}+\rho^{-1}\lambda_k^{(n)}\right) - \rho^{-1}c \\ b \end{bmatrix}.
\end{equation}

    Note that $D$ is a diagonal matrix, because the rows of each matrix $H_k$ are orthonormal, so~\eqref{E:OptCondMinYPrimal} can be solved efficiently, \emph{e.g.}, by block elimination. Moreover, the coefficient matrix is the same at every iteration, so its factorization can be pre-computed and cached before starting the ADMM iterations.

\subsubsection{\textbf{Minimization over $\mathcal{Y}$}} \label{se:MinXblkprimal}
Minimizing~\eqref{E:AugLagrPrimal} over $\mathcal{Y}$ is equivalent to the $p$ independent problems
\begin{equation} 
    \begin{aligned}
    \min_{x_k}  \left\| x_k - H_k x^{(n+1)}  + {\rho}^{-1}\lambda_k^{(n)} \right\|^2 
    \text{subject to } &x_k \in \mathcal{S}_k.
    \end{aligned}
\end{equation}

In terms of the original matrix variables $X_1,\,\ldots,\,X_p$, this amounts to a projection on the PSD cone. 
More precisely, if $\mathbb{P}_k $ denotes the projection onto $\mathbb{S}^{\vert \mathcal{C}_k\vert}_{+}$ we have
\begin{equation} \label{E:XkUpdate}
    x_k^{(n+1)}\!=\!\mathrm{vec}\left\{
    \mathbb{P}_k\!\left[
    \mathrm{vec}^{-1}\!\left( H_k x^{(n+1)} -  {\rho}^{-1}\lambda_k^{(n)}\right)
    \right] \right\}
    .
\end{equation}
Since computing $\mathbb{P}_k$ amounts to an eigenvalue decomposition and each cone $\mathbb{S}^{\vert \mathcal{C}_k\vert}_{+}$ is small by assumption,  we can compute $x_1^{(n+1)},\ldots,x_p^{(n+1)}$ efficiently and in parallel.


\subsubsection{\textbf{Updating the multipliers $\mathcal{Z}$}}
Each multiplier $\lambda_k$, $k=1,\,\ldots,\,p$, is updated with the usual gradient ascent rule,
\begin{equation} \label{E:LambdakUpdate}
    \lambda_k^{(n+1)} = \lambda_k^{(n)} + \rho \left(x_k^{(n+1)} -  H_k x^{(n+1)} \right).
\end{equation}
This computation is cheap, and can be parallelized.

The ADMM algorithm is stopped after the $n$-th iteration if the relative primal/dual error measures $\epsilon_\mathrm{p}$ and $\epsilon_\mathrm{d}$
are smaller than a specified tolerance, $\epsilon_\mathrm{tol}$; see~\cite{boyd2011distributed} for more details on stopping conditions for a generic ADMM algorithm. Algorithm \ref{A:ADMMPrimal} summarizes the the steps to solve a decomposable SDP in standard primal form \eqref{E:PrimalVectorForm}.

\begin{algorithm}[t]
\caption{ADMM for decomposed primal form SDPs}
\label{A:ADMMPrimal}
\begin{algorithmic}[1]
\State \textbf{Input:} $\rho>0$, $\epsilon_\mathrm{tol} >0$, initial guesses $x^{(0)}$, $x_1^{(0)},\ldots,x_p^{(0)}$, $\lambda_1^{(0)},\ldots,\lambda_p^{(0)}$
\State \textbf{Setup:} Chordal decomposition, KKT factorization.
\While{$\max(\epsilon_\mathrm{p}, \epsilon_\mathrm{d}) \geq \epsilon_\mathrm{tol}$}
    \State Compute $x^{(n)}$ with~\eqref{E:OptCondMinYPrimal}.
	\State \textbf{for $k=1,\,\ldots,\,p$:}
		   Compute $x_k^{(n)}$ with~\eqref{E:XkUpdate}.
	\State \textbf{for $k=1,\,\ldots,\,p$:}
		Compute $\lambda_k^{(n)}$ with~\eqref{E:LambdakUpdate}.
	\State Update the residuals $\epsilon_\mathrm{p}, \epsilon_\mathrm{d}$.
\EndWhile
\end{algorithmic}
\end{algorithm}

\section{ADMM for Sparse Dual-form SDPs} \label{se:dSDP}


\subsection{Reformulation of decomposition of the PSD constraint}

Similar to Section~\ref{se:pSDP}, suppose the aggregate sparsity pattern of an SDP in standard dual form~\eqref{E:DualSDP} is described by the chordal graph $\mathcal{G}(\mathcal{V},\mathcal{E})$.
The equality constraint in~\eqref{E:DualSDP} implies that the PSD variable $Z$ has the same sparsity pattern as the aggregate sparsity pattern of the problem data, \emph{i.e.}, $Z\in\mathbb{S}^n_+(\mathcal{E},0)$, so using Theorem~\ref{T:ChordalDecompositionTheorem} we can rewrite~\eqref{E:DualSDP} as
\begin{equation}
\label{E:DualPartiallyDecomposedSDP}
    \begin{aligned}
        \min_{y,Z_1,\ldots,Z_p} \quad & -\langle b,y \rangle\\[-0.5em]
        \text{subject to } \quad
        & \mathcal{A}^*(y)+\sum_{k=1}^p E_k^TZ_kE_k = C,\\
        & Z_k \in \mathbb{S}^{\vert\mathcal{C}_k\vert}_{+}, \quad k=1,\,\ldots,\,p.
    \end{aligned}
\end{equation}
While the original PSD constraint has been replaced by multiple smaller PSD constraints, it is not convenient to apply ADMM to this problem form because the PSD variables $Z_1,\,\ldots,\,Z_k$ in the equality constraint are weighted by the matrices $E_k$.
Instead, we replace $Z_1,\,\ldots,\,Z_k$ in the equality constraint with slack variables $V_1,\ldots,V_p$ such that $Z_k = V_k$, $k = 1,\ldots, p$.
%
Defining $z_k := \text{vec}(Z_k)$ and $v_k := \text{vec}(V_k)$ for all $k=1,\ldots,p$, and using the same vectorized notation as in Section~\ref{se:pSDP} we then reformulate \eqref{E:DualPartiallyDecomposedSDP} in the vectorized form
\begin{equation}\label{E:DualDecomposedSDPVector}
    \begin{aligned}
        \min_{y,z_1,\ldots,z_p,v_1,\ldots,v_p} \quad & -\langle b,y \rangle\\[-0.5em]
        \text{subject to } \quad
        & A^Ty + \sum_{k=1}^p H_k^Tv_k =  c ,\\
        & z_k - v_k = 0, \quad k=1,\,\ldots,\,p,\\
        & z_k \in \mathcal{S}_k,~\qquad\, k=1,\,\ldots,\,p.
    \end{aligned}
\end{equation}

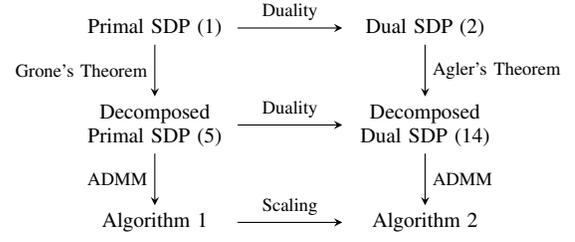
\begin{figure}
    \centering
    \setlength{\abovecaptionskip}{0pt}
    \setlength{\belowcaptionskip}{0em}
	\footnotesize
	\begin{tikzpicture}
	  \matrix (m) [matrix of nodes,
	  		       row sep = 2.5em,	
	  		       column sep = 5em,	
  			       nodes={align=center, text width=2cm}]
  	{
   	Primal SDP~\eqref{E:PrimalSDP} & Dual SDP~\eqref{E:DualSDP} \\
   	Decomposed Primal SDP~\eqref{E:PrimalVectorForm} &
   	Decomposed Dual SDP~\eqref{E:DualDecomposedSDPVector} \\
   	Algorithm 1 & Algorithm 2\\};
	\path[-stealth]
		(m-1-1) edge node [left, align=center]
			{\scriptsize Grone's Theorem}  (m-2-1)
		(m-1-2) edge node [right, align=center]
			{\scriptsize Agler's Theorem}  (m-2-2)
		(m-1-1) edge node [above] {\scriptsize Duality} (m-1-2)
		(m-2-1) edge node [above] {\scriptsize Duality} (m-2-2)
		(m-2-1) edge node [left]  {\scriptsize ADMM}	   (m-3-1)
		(m-2-2) edge node [right] {\scriptsize ADMM}	   (m-3-2)
		(m-3-1) edge node [above] {\scriptsize Scaling} (m-3-2);
	\end{tikzpicture}
    \caption{Duality relationships between primal and dual SDPs, and the decomposed primal and dual SDPs.}
    \label{F:Duality}
\end{figure}

\begin{remark}
    Although we have derived~\eqref{E:DualDecomposedSDPVector} by applying Theorem~\ref{T:ChordalDecompositionTheorem}, \eqref{E:DualDecomposedSDPVector} is exactly the dual of the decomposed primal SDP~\eqref{E:PrimalVectorForm}. Consequently, our analysis provides a decomposition framework for the application of FOMs analogous to that of~\cite{fukuda2001exploiting, kim2011exploiting}  for IPMs. This elegant picture, in which the decomposed SDPs inherit the duality between the original ones by virtue of the duality between Grone's and Agler's theorems, is shown in Fig.~\ref{F:Duality}.
\end{remark}

\subsection{The ADMM algorithm for dual SDPs}

Using indicator functions to move all but the equality constraints $z_k=v_k$, $k=1,\ldots,p$, to the objective gives
\begin{align}
\label{E:DualDecomposedSDPvec}
&\min \,\,
-\langle b,y \rangle + \delta_0\left( c-A^T y-\sum_{k=1}^pH_k^T v_k\right)
+ \sum_{k=1}^{p} \delta_{\mathcal{S}_k}(z_k)
\notag \\
&\text{subject to } \quad
z_k = v_k, \quad k=1,\,\ldots,\,p.
\end{align}
Given a penalty parameter $\rho>0$ and a Lagrange multiplier $\lambda_k$ for each constraint $z_k = v_k$, $k=1,\,\ldots,\,p$, we define
\begin{multline}
\label{E:AugLagrDual}
\mathcal{L}:= -\langle b, y \rangle + \delta_0\left( c - A^T y - \sum_{k=1}^p H_k^T v_k \right)
\\
+ \sum_{k=1}^{p} \left[ \delta_{\mathcal{S}_k}(z_k) + \frac{\rho}{2}\left\| z_k - v_k + \frac{1}{\rho}\lambda_k \right\|^2 \right],
\end{multline}
and group the variables as
$\mathcal{X}:= \{y,v_1,\,\ldots,\,v_p\}$,
$\mathcal{Y}:= \{z_1,\,\ldots,\,z_p\}$, and
$\mathcal{Z}:= \{\lambda_1,\,\ldots,\,\lambda_p\}$.

\subsubsection{\textbf{Minimization over $\mathcal{X}$}}
Minimizing~\eqref{E:AugLagrDual} over block $\mathcal{X}$ is equivalent to the equality-constrained quadratic program
\begin{align}
\label{E:MinYblockDual}
\min_{y,v_1,\ldots,v_p} \quad &-\langle b,y\rangle + \frac{\rho}{2}\sum_{k=0}^{p} \left\| z_k^{(n)} - v_k + \frac{1}{\rho}\lambda_k^{(n)} \right\|^2
\notag\\
\text{subject to } \quad &c - A^T y - \sum_{k=1}^p H_k^Tv_k = 0.
\end{align}
Let $\rho x$ be the multiplier for the equality constraint. After some algebra, the optimality conditions for~\eqref{E:MinYblockDual} can be written as the KKT system
    \begin{equation}
    \label{E:OptCondMinYDual}
    \begin{bmatrix}D & A^T \\ A & 0\end{bmatrix}
    \begin{bmatrix}x \\ y\end{bmatrix} =
    \begin{bmatrix}c - \sum_{k=1}^{p} H_k^T\left( z_k^{(n)}+\rho^{-1}\lambda_k^{(n)}\right) \\ -\rho^{-1}b \end{bmatrix},
    \end{equation}
    plus a set of $p$ uncoupled equations for the variables $v_k$,
    \begin{equation}
    \label{E:vkEqn}
    v_k = z_k^{(n)} +\frac{1}{\rho} \lambda_k^{(n)} + H_k x,  \quad k=1,\,\ldots,\,p.
    \end{equation}
%

The KKT system~\eqref{E:OptCondMinYDual} is the same as~\eqref{E:OptCondMinYPrimal} after rescaling $x\mapsto-x$, $y\mapsto -y$, $c\mapsto \rho^{-1}c$ and $b\mapsto \rho b$, and as in Section~\ref{se:MinYblkprimal} the factors of the coefficient matrix required to solve~\eqref{E:OptCondMinYDual} can be pre-computed and cached. Consequently, updating $\mathcal{X}$ has the same the cost as in Section~\ref{se:MinYblkprimal} plus the cost of~\eqref{E:vkEqn}, which is cheap and can also be parallelized.

\subsubsection{\textbf{Minimization over $\mathcal{Y}$}}
\label{se:MinXblk_dual}

As in Section~\ref{se:MinXblkprimal}, $z_1,\ldots,z_p$ are updated with $p$ independent and efficient projections
\begin{equation}
\label{E:zkUpdate}
z_k^{(n+1)} = \mathrm{vec}\left\{ \mathbb{P}_k \left[ \mathrm{vec}^{-1}\left( v_k^{(n+1)} -  {\rho}^{-1}\lambda_k^{(n)}\right) \right] \right\}.
\end{equation}

\subsubsection{\textbf{Updating the multipliers $\mathcal{Z}$}}
The multipliers $\lambda_k$, $k=1,\,\ldots,\,p$, are updated with the usual gradient ascent rule
\begin{equation}
\label{E:MultUpdateDualSDP}
\lambda_k^{(n+1)} = \lambda_k^{(n)} + \rho\left( z_k^{(n+1)} - v_k^{(n+1)}\right).
\end{equation}


As in Section~\ref{S:ADMMAlgorithmPrimal}, we stop the ADMM algorithm when the relative primal/dual error measures $\epsilon_\mathrm{p}$ and $\epsilon_\mathrm{d}$
%
%
are smaller than a specified tolerance, $\epsilon_\mathrm{tol}$. Algorithm~\ref{A:ADMMDual} summarizes the full ADMM algorithm for sparse dual-standard-form SDPs.


\begin{algorithm}[t]
    \caption{ADMM for decomposed dual form SDPs}
    \label{A:ADMMDual}
    \begin{algorithmic}[1]
\State \textbf{Input:} $\rho>0$, $\epsilon_\mathrm{tol} >0$, initial guesses $y^{(0)}$, $z_1^{(0)},\ldots,z_p^{(0)}$, $\lambda_1^{(0)},\ldots,\lambda_p^{(0)}$
\State \textbf{Setup:} Chordal decomposition, KKT factorization.
\While{$\max(\epsilon_\mathrm{p}, \epsilon_\mathrm{d}) \geq \epsilon_\mathrm{tol}$}
    	\State \textbf{for $k=1,\,\ldots,\,p$:}
    			Compute $z_k^{(n)}$ with~\eqref{E:zkUpdate}.
    	\State Compute $y^{(n)},x$ with~\eqref{E:MinYblockDual}.
    	\State \textbf{for $k=1,\,\ldots,\,p$:}
    			Compute $v_k^{(n)}$ with~\eqref{E:vkEqn}
    	\State Compute $\lambda_k^{(n)}$ with~\eqref{E:DualMultUpdate} (no cost).
    	\State Update the residuals $\epsilon_\mathrm{p}$ and $\epsilon_\mathrm{d}$.
\EndWhile
\end{algorithmic}
\end{algorithm}
%


\begin{remark}
    The computational cost of~\eqref{E:vkEqn} is the same as~\eqref{E:LambdakUpdate}, 
    so the ADMM iterations for the decomposed dual-standard-form SDP~\eqref{E:DualDecomposedSDPVector} have the same cost as those for the decomposed primal-standard-form SDP~\eqref{E:PrimalVectorForm}, plus the cost of~\eqref{E:MultUpdateDualSDP}. However, if one minimizes~\eqref{E:AugLagrDual} over $\mathcal{Y}$ \emph{before} minimizing it over $\mathcal{X}$, substituting \eqref{E:vkEqn}  into~\eqref{E:MultUpdateDualSDP} gives
    \begin{equation} \label{E:DualMultUpdate}
        \lambda_k^{(n+1)} =  \rho H_k x^{(n+1)}, \quad k=1,\,\ldots,\,p.
    \end{equation}
Since $H_1 x,\,\ldots,\,H_px$ have already been computed to update $v_1,\,\ldots,\,v_p$, updating $\lambda_1,\,\ldots,\,\lambda_p$ requires only a scaling operation. Consequently, the ADMM algorithms for the primal- and dual-standard-form SDPs can be considered as scaled versions to each other, with the same leading-order computational cost at each iteration.
\end{remark}

\section{Numerical Simulations} \label{se:simulation}

We have implemented our techniques in CDCS (Cone Decomposition Conic Solver)~\cite{CDCS}, an open-source MATLAB solver. CDCS supports cartesian products of the following cones: $\mathbb R^n$, non-negative orthant, second-order cone, and the PSD cone. Currently, only chordal decomposition techniques for semidefinite cones are implemented, while the other cone types are not decomposed. 
Our codes can be downloaded from
%
{\small
\url{https://github.com/OxfordControl/CDCS}}
\normalsize.

We tested CDCS on four sparse large-scale ($n \geq 1000, m \geq 1000$) problems in SDPLIB \cite{borchers1999sdplib}, as well as on randomly generated SDPs with block-arrow sparse pattern, used as a benchmark in~\cite{sun2014decomposition}. The performance is compared to that of the IPM solver SeDuMi \cite{sturm1999using} and of the first-order solver SCS~\cite{scs} on both the full SDPs (without decomposition) and the SDPs decomposed by SparseCoLO~\cite{fujisawa2009user}.

The comparison has two purposes: 1) SeDuMi computes accurate optimal points, which can be used to assess the quality of the solution computed by CDCS; 2) SCS is a high-performance first-order solver for general conic programs, so we can assess the advantages of chordal decomposition. SeDuMi should not be compared to the other solvers on CPU time, because the latter only aim to achieve moderate accuracy. 
In the experiments reported below, the termination tolerance for CDCS and SCS was set to $\epsilon_\mathrm{tol} = 10^{-3}$, with a maximum of 2000 iterations. All experiments were carried out on a PC with an Intel(R) Core(TM) i7 CPU, 2.8 GHz processor and 8GB of RAM.

\subsection{SDPs with block-arrow pattern}

\begin{figure}[t]
\centering
\setlength{\abovecaptionskip}{0pt}
\setlength{\belowcaptionskip}{0em}
\tikzset{decorate sep/.style 2 args=
{decorate,decoration={shape backgrounds,shape=circle,shape size=#1,shape sep=#2}}}
\begin{tikzpicture}
	 \draw (0,0) rectangle  (2.5,2.5);
	 \draw[fill=black!20] (0,2.5) rectangle  (0.5,2);
	 \draw[fill=black!20] (0.5,2) rectangle  (1,1.5);
	 \draw[fill=black!20] (0,0)--(2.5,0)--(2.5,2.5)--
	 				      (2.15,2.5)--(2.15,0.35)--(0,0.35)--(0,0);
	 \draw[decorate sep={0.5mm}{2mm},fill] (1.2,1.3)--
	 					   node[left] {\scriptsize $l$ blocks $ $ }(2.0,0.5);
	 \draw[<->] (0,2.6) -- node[above] {\scriptsize $d$} (0.5,2.6) ;
	 \draw[<->] (-0.1,2) -- node[left] {\scriptsize $d$} (-0.1,2.5) ;
	 \draw[<->] (2.15,2.6) -- node[above] {\scriptsize $h$} (2.5,2.6) ;
	 \draw[<->] (-0.1,0) -- node[left] {\scriptsize $h$} (-0.1,0.35) ;
\end{tikzpicture}
\caption{Block-arrow sparsity pattern: the number of blocks, $l$; block size, $d$; the size of the arrow head, $h$.}
\label{F:BlockArrowSDP}
\end{figure}
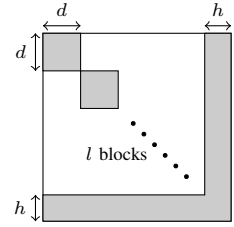

SDPs with the block-arrow sparsity pattern shown in Fig.~\ref{F:BlockArrowSDP}---which is chordal---are used as a benchmark in \cite{sun2014decomposition}. The SDP parameters are: the number of blocks, $l$; the block size, $d$; the size of the arrow head, $h$; the number of constraints, $m$. Here, we consider the following cases:
1) Fix $l=40$, $d=10$, $h=20$, vary $m$;
2) Fix $m=1000$, $d=10$, $h=20$, vary $l$;
3) Fix $l=40$, $h=10$, $m=1000$, vary $d$.

\begin{figure*}[t]
    \centering
    \setlength{\abovecaptionskip}{0pt}
    \setlength{\belowcaptionskip}{0em}
	\includegraphics[scale=1]{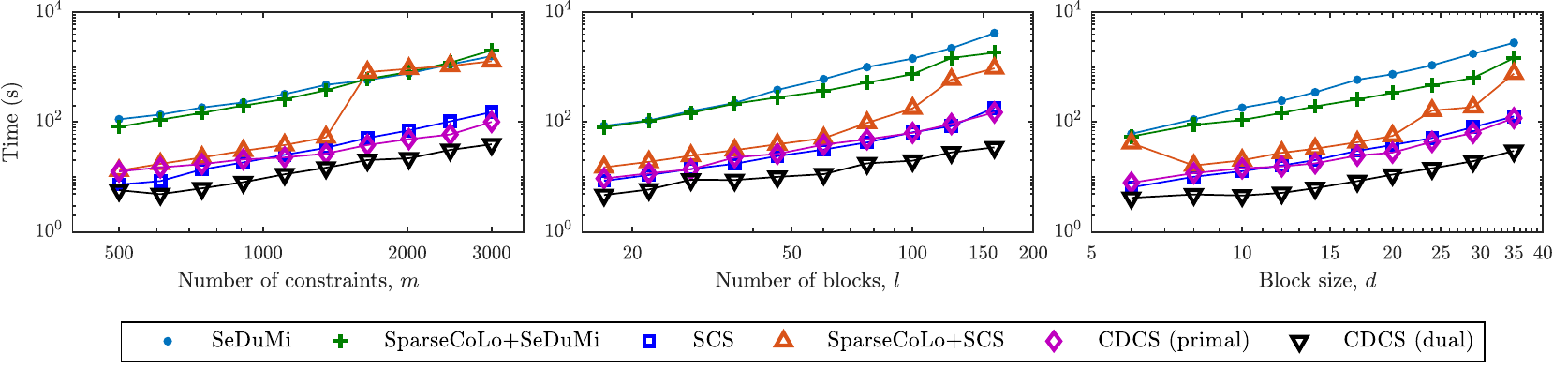}
    \caption{CPU time for SDPs with block-arrow patterns. Left to right: varying number of constraints; varying number of blocks; varying block size.}
    \label{F:BlcokArrowResult}
\end{figure*}

The CPU times for different solvers, averaged over five random problem instances, are shown in Fig.~\ref{F:BlcokArrowResult}. CDCS is approximately 10 times faster than SeDuMi and the combination SparseCoLO+SeDuMi, our Algorithm~\ref{A:ADMMDual} being the fastest. Besides, the optimal value from CDCS was always within 0.02\% of the accurate value from SeDuMi.

\subsection{Sparse SDPs from SDPLIB}
We consider two max-cut problems (maxG32 and maxG51), a Lov\'asz theta problem (thetaG51), and a box-constrained quadratic problem (qpG51) from SDPLIB \cite{borchers1999sdplib}. Table~\ref{T:SparseStatistic} reports the dimensions and chordal decomposition details of these large, sparse SDPs. 
Table \ref{T:ResultsSDPLIB} summarizes our numerical results; maxG51, thetaG51 and qpG51 could not be solved by the combination SparseCoLO+SeDuMi due to memory overflow. 
For all four problems, CDCS (primal and dual) is faster than SeDuMi and can give speedups compared to SCS and SparseCoLO+SCS. 
We remark that the stopping objective value from CDCS is within 2\% of the optimal value returned by SeDuMi (which is highly accurate, and can be considered exact) in all four cases, and within 0.08\% for the max-cut problems maxG32 and maxG51 --- a negligible difference in applications. 
\begin{table}[t]
    \centering
    \setlength{\abovecaptionskip}{0pt}
    \setlength{\belowcaptionskip}{0em}
    \renewcommand\arraystretch{0.9}
    \caption{Problem statistics for SDPLIB problems}
    \label{T:SparseStatistic}
    \begin{tabular}{r c c c c}
        \hline \toprule[1pt] 
            & maxG32 & maxG51  & thetaG51 & qpG51 \\
        \hline\\[-0.75em]
        Affine constraints, $m$         & 2000 & 1000 & 6910 & 1000 \\
        Original cone size, $n$         & 2000 & 1000 & 1001 & 2000 \\
        Number of cliques, $p$  & 1499 & 674 & 674 & 1675   \\
        Maximum clique size  & 60 & 322 & 323 & 304   \\
        Minimum clique size  & 5 & 6 & 7 & 1   \\
        \bottomrule[1pt]
        \end{tabular}
\end{table}

\begin{table*}
\centering
\renewcommand\arraystretch{0.75}
\setlength{\abovecaptionskip}{0pt}
\setlength{\belowcaptionskip}{0em}
\caption{Results for the problem instances in SDPLIB}
\label{T:ResultsSDPLIB}
\begin{tabular}{c r c c c c c c}
\toprule[1pt]
& & SeDuMi
&\begin{tabular}[x]{@{}c@{}}SparseCoLO+\\SeDuMi\end{tabular}
& SCS
&\begin{tabular}[x]{@{}c@{}}SparseCoLO+\\SCS\end{tabular}
&\begin{tabular}[x]{@{}c@{}}CDCS\\(primal)\end{tabular}
&\begin{tabular}[x]{@{}c@{}}CDCS\\(dual)\end{tabular}\\
\midrule
\multirow{4}{*}{maxG32} &
Total time (s) 			 & 974.6 & 355.2 & 2.553 $\times 10^3$ & 65.1  & 88.6 & 53.1\\
& Pre-processing time (s)&  0    & 3.18 & 0.43  & 3.24 & 21.2 & 21.4\\
& Objective value        & 1.568$\times 10^3$ & 1.568$\times 10^3$ & 1.568$\times 10^3$ & 1.566$\times 10^3$  & 1.569$\times 10^3$ & 1.568$\times 10^3$  \\
& Iterations             & 14 & 15 & 2000 & 960 & 238  & 127 \\
\midrule
\multirow{4}{*}{maxG51}  &
Total time (s)           & 134.5 & -- & 87.9 &  1.201$\times 10^3$ & 110.9 & 75.9\\
& Pre-processing time (s)&  0   & -- & 0.11  &  2.87 & 3.30 & 3.20 \\
& Objective value        & 4.006$\times 10^3$ & -- & 4.006$\times 10^3$ & 3.977$\times 10^3$  & 4.005$\times 10^3$  & 4.006$\times 10^3$ \\
& Iterations             & 16 & -- & 540 & 2000 & 235  & 157 \\
\midrule
\multirow{4}{*}{thetaG51}&
Total time (s)           & 2.218 $\times 10^3$ & -- & 424.2 & 1.346 $\times 10^3$    & 471.2  & 735.1 \\
& Pre-processing time (s)&  0   & -- & 0.30  & 5.30 & 25.1 & 25.0\\
& Objective value        & 349 & -- & 350.6 & 341.3  & 354.5 & 355.9 \\
& Iterations             & 20 & -- & 2000 & 2000   & 394  & 646 \\
\midrule
\multirow{4}{*}{qpG51} &
Total time (s)           & 1.407$\times 10^3$ & -- & 2.330$\times 10^3$ & 985.8   & 727.1  & 606.2 \\
& Pre-processing time (s)&    0    & -- & 0.47  & 190.2  & 12.3  & 12.3 \\
& Objective value        & 1.182$\times 10^3$ & -- & 1.288$\times 10^3$ & 1.174$\times 10^3$   & 1.195$\times 10^3$  & 1.194$\times 10^3$  \\
& Iterations             & 22 & -- & 2000 & 2000  & 1287   & 1048   \\
\bottomrule[1pt]
\end{tabular}
\vspace*{-1.5em}
\end{table*}

\section{Conclusion} \label{se:conclusion}

We proposed a conversion framework for SDPs characterized by chordal sparsity suitable for the application of FOMs, analogous to the conversion techniques for IPMs of~\cite{fukuda2001exploiting, kim2011exploiting}. We also developed efficient ADMM algorithms for sparse SDPs in primal or dual standard form, which are implemented in the solver CDCS.
Numerical experiments on SDPs with block-arrow sparsity patterns and on large sparse problems in SDPLIB show that our methods can provide speedups compared to both IPM solvers such as SeDuMi~\cite{sturm1999using}---even when the chordal sparsity is exploited using SparseCoLO~\cite{fujisawa2009user}---and the state-of-the-art first-order solver SCS~\cite{scs}.
Exploiting chordal sparsity in a first-order HSDE formulation similar to that of~\cite{ODonoghue2016} would be desirable in the future to be able to detect infeasibility. 


\balance
\bibliographystyle{IEEEtran}
\bibliography{Reference}
 \end{document}